\documentclass[a4paper,reqno]{amsart}

\usepackage[utf8]{inputenc}
\usepackage[T1]{fontenc}
\usepackage[a4paper]{geometry}

\usepackage{amssymb}

\usepackage{stmaryrd}
\usepackage[cal=boondoxo,scr=pxtx]{mathalfa}
\usepackage{amsthm}
\theoremstyle{plain}
\newtheorem{theorem}{Theorem}
\newtheorem{property}{Property}

\newtheorem{lemme}[theorem]{Lemma}
\newtheorem{corollaire}[theorem]{Corollary}

\theoremstyle{definition}
\newtheorem{definition}{Definition}

\theoremstyle{remark}
\newtheorem*{remark}{Remark}

\usepackage[colorlinks,hyperindex,unicode]{hyperref}  %
\hypersetup{
    colorlinks,
    citecolor=blue,
    linkcolor=blue,
    urlcolor=blue,
}

\usepackage{orcidlink}
\usepackage{booktabs}

\newcommand{\llbra}{\llbracket}
\newcommand{\rrket}{\rrbracket}
\newcommand{\sgn}{\mathrm{sgn}}
\newcommand{\Lie}{\ensuremath{\mathscr{L}}}
\newcommand{\sss}[1]{\scriptscriptstyle{#1}}
\newcommand{\dr}{\partial}
\renewcommand{\S}{\ensuremath{\mathscr{S}}}
\renewcommand{\div}{\ensuremath{\mathrm{div}}}
\newcommand{\sS}{{\sss S}}
\newcommand{\sA}{{\sss A}}
\newcommand{\sB}{{\sss B}}
\newcommand{\sC}{{\sss C}}
\newcommand{\sn}{{\mathfrak{n}}}
\newcommand{\eP}{\epsilon_\perp}
\newcommand{\Op}{\mathrm{Op}}
\newcommand{\g}{{\scriptstyle{\mathsf{g}}}}  %
\newcommand{\grad}{\g}
\renewcommand{\tilde}{\widetilde}
\renewcommand{\Im}{\operatorname{Im}}
\def\L{\Lie}
\def\setZ{\mathbb{Z}}
\def\R{\mathscr{R}}
\def\Id{\mathrm{{Id}}}
\def\Tr{\mathrm{Tr}}

\usepackage{amsaddr}

\begin{document}
    \title[Submanifolds, restriction of differential operators 
    ...]{Submanifolds, restriction of differential operators and 
    Schouten-Nijenhuis bracket}
    \date{\today}
    \author[J.~Queva]{J.~Queva \orcidlink{0000-0001-8280-1925}}
    \address{Universit\'e de Corse -- CNRS UMR 6134 SPE, Campus Grimaldi BP 52, 
    20250 Corte, France.}
    \email{queva@univ-corse.fr}
    \author[J.~Renaud]{J.~Renaud \orcidlink{0000-0002-3284-5003}}
    \address{APC-Astroparticule et Cosmologie (UMR-CNRS 7164), 
        Batiment Condorcet, 10 rue Alice Domon et L\'eonie Duquet, F-75205 
        Paris Cedex 13, France.}
    \email{jrenaud@apc.in2p3.fr}

    \begin{abstract}
        In this article, we present new formulas for restricting the Hodge 
        operator, the codifferential and the Hodge Laplacian to a submanifold.
        These formulas are concise, thanks to the use of tools for differential 
        forms (Schouten-Nijenhuis bracket and generalized interior product). 
        They also facilitate explicit calculations for practical applications.
    \end{abstract}
    
    \maketitle
    
    \tableofcontents

    \section{Introduction}
    
    Here, we present restriction formulas to a submanifold for the Hodge 
    operator [Eq.~\eqref{Eq. hodge}], the codifferential [Eqs.~\eqref{Eq. 
    delta1} and \eqref{Eq. delta}] and the Hodge Laplacian [Eqs.~\eqref{Eq. 
    square} and \eqref{Eq. square1}] acting on $p$-forms.
    These formulas are concise and facilitate practical computations for 
    low-degree differential forms [Eqs.~\eqref{eq:First}--\eqref{eq:Last}].

    Recently, in Ref.~\cite{Chan-Czubak},  the question of restricting the 
    Hodge Laplacian ($\Delta=-(d\delta+\delta d)$) to a submanifold has been 
    addressed, limited to the case of the $1$-form field.
    This work  relies on decomposing the Levi-Civita connection $\nabla$ into 
    tangential and normal components, a decomposition that essentially dates 
    back to Gauss (theorema egregium, see, e.g., 
    Refs.~\cite{Jost,Frankel}).
    This approach yields the restriction of the rough Laplacian 
    ($\nabla^*\nabla=\Tr(\nabla^2)$) on 1-forms, and the Hodge Laplacian is 
    then obtained using the Weitzenböck formula.
    This question is drawing renewed interest in fluid mechanics, see 
    Ref.~\cite{Chan-Czubak-Fuster}.

    Another approach  is to restrict the Hodge operator first, then the 
    codifferential, and finally the Hodge Laplacian.
    This method has received little attention.
    Here, we use it to obtain new compact formulas for the restriction, 
    regarding differential forms of arbitrary degree.
    This relies on two key tools: generalized interior and exterior products on 
    forms (the $i$ and $j$ operators defined below), and the correspondence 
    between two generalizations of the Lie derivative (the Schouten-Nijenhuis 
    bracket 
    \cite{Schouten:1953,Nijenhuis:1955A,Nijenhuis:1955B,Marle97,Huguet:2024mza} 
    and the Michor-Tulczyjew operator \cite{Michor1987,Tulczyjew:1974}).
    The formulas thus obtained involve various objects from the geometry of 
    differential forms on a pseudo-Riemannian manifold.
    It is noteworthy that these compact formulas allow for explicit 
    calculations for forms of low degree.

    The paper is organized as follows.
    Section~\ref{Sec. Prelim} presents the necessary tools for working with 
    differential forms.
    Section~\ref{Sec. Context} sets the geometric context.
    Section~\ref{Sec. Restr.} presents the main results of the paper: the 
    restriction of the Hodge operator, the codifferential and the Hodge 
    Laplacian.
    Section~\ref{Sec. transv} is concerned with a specific kind of differential 
    forms, called transverse, for which the restriction formulas simplify 
    significantly.
    Finally, in Sec.~\ref{Sec. small a}, we combine all the preceding results 
    to compute explicit formulas for low-degree differential forms.

    \section{Required tools for working with differential forms}\label{Sec. 
    Prelim}

    In this section, we consider an oriented pseudo-Riemannian differential 
    manifold $(E,g)$ of finite dimension.
    We denote by $\Omega(E)$ the set of alternating forms on $E$ and 
    $\Omega^a(E)$ the set of alternating forms of degree $a$.
    We denote by $\hat\eta$ the operator defined by 
    $\hat\eta(\alpha)=(-1)^a\alpha$ when $\alpha\in\Omega^a(E)$.
    The metric implies a natural isomorphism 
    $\flat: \Lambda T(E)\to\Lambda T^*(E)=\Omega(E)$ and its inverse  
    $\sharp=\flat^{-1}$.
    
    An operator $X$ on the space of alternating forms $\Omega(E)$ is said to be 
    of degree $x\in\setZ$ if $X\Omega^a(E)\subset\Omega^{a+x}(E)$.
    We denote by $\Op^x$ the set of operators of degree $x$.
    We also use the graded commutator defined by:
    \begin{equation}
        \label{eq:Graded_Bracket_1}
        [X,Y]^{\grad}=XY-(-1)^{xy}YX,
    \end{equation} 
    where $X\in \Op^x$ and $Y\in\Op^y$.
    It satisfies the graded identities:
    \begin{align}
        &[X,Y]^\g = -(-1)^{xy}[Y,X]^\g,\\
        &[X,YZ]^\g = [X,Y]^\g Z + (-1)^{xy}Y[X,Z]^\g\label{eq:Graded_Bracket_3},
    \end{align}
    where $X\in \Op^x$, $Y\in\Op^y$ and $Z\in\Op^z$.
    
    \subsection{A Conjugation Operator}
    
    We now introduce a conjugation operator defined relative to the Hodge 
    operator, which was used implicitly in Ref.~\cite{Huguet:2024mza} and which 
    we make explicit from this point onward.
    \begin{definition}
        For $X\in\Op^x$, we set
        \begin{equation*}
            \R (X)= *^{-1}X*\hat\eta^x,
        \end{equation*}
        where $*$ is the Hodge operator. That is, for $\alpha\in\Omega^a$
        \begin{equation*}
            \R (X)\alpha=(-1)^{ax} *^{-1}X*\alpha. 
        \end{equation*}
    \end{definition}
    This operator fulfills the following properties.
    \begin{property}\label{prop:R}
        For $X\in \Op^x$ and $Y\in\Op^y$, we have
        \begin{align}
            &\R(X)\in\Op^{-x},\nonumber\\ %
            &\R(XY)=(-1)^{xy}\R(X)\R(Y),\nonumber\\ %
            &\R^2(X)=(-1)^xX,\nonumber\\ %
            &\R^{-1}(X)=*^{-1}X*(-\hat\eta)^x=(-1)^x\R(X),\nonumber\\ %
            &\R([X,Y]^\g) = -[\R(Y),\R(X)]^\g.\label{eq.R_Commutateur}
        \end{align}
    \end{property}
    \begin{proof}
        Direct calculation from the definition.
        Note that the $\hat{\eta}^x$ term in the definition ensures that the 
        formulas are dimension-independent.
    \end{proof}

    As a direct application, we have
    \begin{equation*}
        \delta=\R(d),\quad 
        \R(\delta)=-d,\quad
        \R(\Delta)=\Delta,\quad
        \R(\nabla_v) = \nabla_v,
    \end{equation*}
    where the first identity relies on the fact that $\R(d)$ coincides with the 
    definition of $\delta$ since $d\in\Op^1$, the following two are 
    consequences of Prop.~\ref{prop:R}.
    The last identity involves the Levi-Civita connection $\nabla$ which 
    belongs to $\Op^0$ and commutes with the Hodge operator.

    \subsection{Exterior and Interior Products on Forms}
    Following Ref.~\cite{Huguet:2024mza}, the exterior product $j$ and the 
    interior product $i$ are defined by
    \begin{equation*}
        j^\alpha\beta=\alpha\wedge\beta,\quad
        i^\alpha=\R^{-1}(j^\alpha),\quad
        \forall\alpha,\beta\in\Omega(E),
    \end{equation*}
    with $j^\alpha\in\Op^a$ and $i^\alpha\in\Op^{-a}$ for $\alpha\in\Omega^a$.
    The properties of these operators are detailed in 
    Ref.~\cite{Huguet:2024mza}.
    Here we will mainly use that $i^{\alpha\wedge\beta}=i^\beta i^\alpha$, 
    which is immediately deduced from the definition.
    These definitions generalize to alternating multi-vectors by
    \begin{equation*}
        i_u=i^{\flat u} \ \mbox{ and }\ j_u=j^{\flat u},\quad
        \forall u\in\Lambda T(E).
    \end{equation*}
    In the case where $u$ is of degree one, we recover the usual definitions of 
    interior and exterior products (see, e.g., Ref.~\cite{Fecko} Ch.~5) and we 
    have, if $\lambda\in\Omega^1(E)$  and $u\in T(E)$,
    \begin{equation*}
        i_u j^\lambda+j^\lambda i_u
        =\langle\lambda,u\rangle,%
    \end{equation*}
    where $\langle\ ,\ \rangle$ is the natural pairing between forms and 
    vectors.
    For a given frame $e_a$ and its coframe $e^a$ we denote $i_a \equiv 
    i_{e_a}$, $i^a \equiv i^{e^a}$, $j_a \equiv j_{e_a}$, and $j^a \equiv 
    j^{e^a}$.
    
    \subsection{Two Generalizations of the Lie Derivative and Their 
    Correspondence}
    
    The Schouten-Nijenhuis bracket extends the Lie derivative to alternating 
    multivectors (see 
    Refs.~\cite{Schouten:1953,Nijenhuis:1955A,Nijenhuis:1955B} and 
    Ref.~\cite{Marle97} for an axiomatic presentation).
    Moreover, the identification between  $ \Lambda T(E)$ and $\Lambda T^*(E)$ 
    (through the musical applications $\flat$ and $\sharp$) enables the 
    transport of this bracket to forms.
    We denote this transported bracket by $\llbracket\ , \ \rrbracket$, that is 
    to say 
    $\llbracket\alpha,\beta\rrbracket=\flat[\sharp\alpha,\sharp\beta]_{\rm 
    SN}$, where $[\ ,\ ]_{\rm{SN}}$ is the usual Schouten-Nijenhuis bracket 
    (see Ref.~\cite{Huguet:2024mza} for a detailed presentation and a wealth of 
    identities involving it).
    This bracket measures the failure of the codifferential to be a derivation 
    with respect to the exterior product $\wedge$, since we have (see 
    Ref.~\cite{Huguet:2024mza} again and references therein), for  
    $\alpha\in\Omega^a$ and $\beta\in\Omega^b$,
    \begin{equation}\label{delta_j}
        \delta(\alpha\wedge\beta)= 
        \delta(\alpha)\wedge\beta+(-1)^a\alpha\wedge\delta(\beta)+(-1)^a\llbracket\alpha,\beta\rrbracket.
    \end{equation}
    Note that, unlike the historical Schouten-Nijenhuis bracket, this bracket 
    depends on the metric since the operators $\sharp$ and $\flat$ do.
    In the following, we will often write
    \begin{equation*}
        \S^\alpha\beta=\llbracket\alpha,\beta\rrbracket,
    \end{equation*}
    with $\S^\alpha\in \Op^{a-1}$, which allows us to rewrite 
    Eq.~\eqref{delta_j} in the form
    \begin{equation}
        [\delta, j^\alpha]^{\grad}
        = j^{\delta\alpha}+(-1)^a\S^\alpha,\
        \forall\alpha\in\Omega^a(E).\label{deltaj}
    \end{equation}
    Though unnecessary here, we record the following lemma for possible use 
    elsewhere.
    \begin{lemme}\label{Lem:Lemma1}
        For $\lambda\in\Omega^1$, $\nabla$ the Levi-Civita connection and 
        $\L_v$ the Lie derivative along $v$, we have
        \begin{equation*}%
            \S^\lambda\alpha
            =-[(\nabla\lambda)_{ml}
            + (\nabla\lambda)_{lm}] j^m i^l\alpha+\Lie_{\sharp \lambda}\alpha.
        \end{equation*}
    \end{lemme}
    \begin{proof}
        This follows directly from Eq.~(9) in \cite{Huguet:2024mza}, $(\L_w 
        g)(u,v) = g(\nabla_u w, v) + g(\nabla_v w, u)$, and the identity 
        $\Lie_v(\sharp g) = -\sharp\Lie_v g$.
    \end{proof}

    Another generalization of the Lie derivative is given by the 
    Michor-Tulczyjew operator \cite{Michor1987,Tulczyjew:1974} which is 
    defined, for $\alpha\in\Omega^a$, by
    \begin{equation*}
        \Theta^\alpha=i^\alpha d-(-1)^a d i^\alpha=[i^\alpha, d]^{\grad}\in 
        \Op^{1-a},
    \end{equation*}
    that extends to multi-vectors of degree $a$ by $\Theta_u=i_u d-(-1)^adi_u$, 
    and coincides with the Lie derivative if $a=1$.
    We then have the following property, which will play a decisive role in the 
    restriction of the codifferential.
    \begin{property}\label{R_Theta}
        The Schouten-Nijenhuis and Michor-Tulczyjew operators correspond to 
        each other under the transformation $\R$ according to
        \begin{align}
            &\R(\Theta^\alpha)
            = -(j^{\delta\alpha}+(-1)^a\S^\alpha),\quad
            \forall\alpha\in\Omega^a,\label{Eq. R_Theta}\\
            &\R(\S^\alpha)=i^{\delta\alpha}+\Theta^\alpha.\label{Eq. R_Schouten}
        \end{align}
    \end{property}
    \begin{remark}
        The special case $a=1$ yields the following identity:
        \begin{equation*}
            \R(\Lie_v)=\S^{\flat v}+\div(v),
            \quad\forall v\in T(E),
        \end{equation*}
        with $\L_v$ the Lie derivative along $v$ and $\div(v)=-\delta(\flat 
        v)$  the divergence of the vector field $v$.
    \end{remark}
    \begin{proof}
        For Eq.~\eqref{Eq. R_Theta} a direct computation involving the previous 
        identities and Eq.~\eqref{deltaj} in particular.
        Equation~\eqref{Eq. R_Schouten} is obtained by applying $\R$ to 
        Eq.~\eqref{Eq. R_Theta}.
    \end{proof}
    
    The detailed examples above are summarized in Table~\ref{tab:exemple}.

    \begin{table}
        \centering
        \begingroup
        \renewcommand{\arraystretch}{1.4}
        \begin{tabular}{l @{\qquad} l}
            \toprule
            $i^\alpha=\R^{-1}(j^\alpha)=(-1)^a\R(j^\alpha)$  
            &$\R(i^\alpha)=j^\alpha$  \\
            $\delta=\R(d)$ & $\R(\delta)=-d$ \\
            $\R(\Delta)=\Delta$ & $\R(\nabla_u)=\nabla_u$  \\
            $\R(\Theta^\alpha)=-(j^{\delta\alpha}+(-1)^a\S^\alpha)$&$\R(\S^\alpha)=\Theta^\alpha+i^{\delta\alpha}$\\
            \bottomrule
        \end{tabular}
        \endgroup
        \caption{\label{tab:exemple} 
            Recording of various formulas involving the conjugation operator 
            $\R$. 
        }
    \end{table}
    
    Furthermore, from well-known or easily established commutators, we can 
    obtain, by applying Eq.~\eqref{eq.R_Commutateur}, less obvious commutators 
    that we use crucially in Sec.~\ref{Sec. transv} (see a summary in 
    Table~\ref{tab:commutateurs}).
    For instance, since $[d,\Theta^\alpha]^{\grad}=0$, we recover, using 
    Eq.~(\ref{deltaj}) and $\delta^2=0$, that 
    $[\delta,\S^\alpha]^{\grad}=\S^{\delta\alpha}$, in other words
    \begin{equation*}
        \delta(\llbracket\alpha,\beta\rrbracket)=\llbracket\delta\alpha,\beta\rrbracket-(-1)^a\llbracket\alpha,\delta\beta\rrbracket,
    \end{equation*}
    and $\delta$ is, ultimately, a derivation, but with respect to the 
    Schouten-Nijenhuis bracket.
    
    \begin{table}\label{tab gradcom}
        \centering
        \begingroup
        \renewcommand{\arraystretch}{1.4}
        \begin{tabular}{l @{\qquad} l}
            \toprule
            $[j^\alpha,j^\beta]^\g = 0$
            & $[i^\alpha, i^\beta]^\g = 0$\\
            $[d,j^\alpha]^\g = j^{d\alpha}$ 
            & $[i^\alpha,\delta]^\g = i^{d\alpha}$ \\
            $[\delta,j^\alpha]^\g = j^{\delta\alpha} + (-1)^a\S^\alpha$
            & $[i^\alpha, d]^\g = \Theta^\alpha$\\
            $[\nabla_u,j^\alpha]^\grad=j^{\nabla_u\alpha}$
            & $[i^\alpha,\nabla_u]^\grad=-i^{\nabla_u\alpha}$\\
            $[\Theta^\alpha,d]^\g = 0$
            & $[\delta,\S^\alpha]^\g = \S^{\delta\alpha}$ \\
            $[\Theta^\alpha,i^\beta]^\g = - i^{\llbra\beta,\alpha\rrket}$
            & $[\S^\alpha,j^\beta]^\g = j^{\llbra\alpha,\beta\rrket }$\\
            $[\Theta^\alpha,\Theta^\beta]^\g
            = -\Theta^{\llbra \beta,\alpha\rrket}$
            &  $[\S^\alpha,\S^\beta]^\g = \S^{\llbra\alpha,\beta\rrket}$ \\
            \bottomrule
        \end{tabular}
        \endgroup
        \caption{\label{tab:commutateurs} 
            Recording of various commutators related by the conjugation $\R$.
        }
    \end{table}

    \section{Geometric context of this article}
    \label{Sec. Context}
    
    \emph{In all that follows, we will only consider local properties, and all 
    objects will be defined only locally.}
    
    We assume that $S$ is an oriented submanifold of dimension $n$ embedded in 
    an oriented pseudo-Riemannian manifold $(M, g)$ of dimension $n+k$ via the 
    canonical injection $s$.
    We denote $g_\sS$ the induced metric on $S$ and we assume that $T(S)$ is 
    non-isotropic, which implies that $S$ is also a pseudo-Riemannian manifold.
    We denote $S^\perp$ the normal bundle whose base is $S$ and whose fiber at 
    the point $y \in S$ is $T_y(S)^\perp = \{u\in T_y(M)\ |\  g(u,v)=0,\ 
    \forall v\in T_y(S)\}$.
    The bundle $S^\perp$ inherits a metric $g^\perp$, an orientation, and a 
    volume form $\omega^\sn$, and
    if $e^{n+1}, \dots, e^{n+k}$ is a direct orthonormal co-basis of $S^\perp$, 
    then $\omega^\sn = e^{n+1} \wedge \cdots \wedge e^{n+k}$.
    Additionally, we denote $\omega_\sn$  the $k$-vector volume of $S^\perp$: 
    $\omega_\sn =\sgn(g^\perp)\sharp\omega^\sn= e_{n+1} \wedge \cdots \wedge 
    e_{n+k}$, where, as usual, the $e_a$ are the (pre-)dual basis of the $e^a$.

    We work under the following hypotheses which we will call (H) in the rest 
    of the article:
    \begin{itemize}
        \item[\hypertarget{H}{(H)}] Let $S\subset M$ be two manifolds of 
        respective dimensions $n$ and $n+k$ and $s : S\to M$ be the canonical 
        injection. $M$ and $S$ are respectively equipped with the metrics $g$ 
        and $g_\sS=s^*g$ for which they are both oriented pseudo-Riemannian 
        manifolds.
        Let $\omega^\sn$ be the volume form of the normal bundle $S^\perp$ and 
        $\omega_\sn=\sgn(g^\perp)\sharp\omega^\sn$.\label{(H)}
    \end{itemize}

    We call the restriction (to $S$) of $\alpha \in \Omega(M)$ the pull-back of 
    $\alpha$ by $s$, denoted $s^*\alpha$.
    The end goal of this article is to compare the restriction of the Laplacian 
    with the Laplacian of the restriction: $s^*(\Delta\alpha) = \Delta_\sS 
    (s^*\alpha) + T$, where $T$ is the additional term to be calculated, and, 
    as everywhere in this article, operators related to the submanifold $S$ are 
    indexed by $S$, while operators without indices relate to the ambient 
    manifold $M$.
    Since $\Delta=-(d\delta+\delta d)$ and $\delta=*^{-1}d*\hat\eta$, we will 
    proceed by evaluating $s^**$, $s^*\delta$ and, finally, $s^*\Delta$, given 
    that $s^*d$ is already known: 
    \begin{equation*}
        s^*d = d_\sS s^*,
    \end{equation*}
    see, e.g., Eq.~6.2.12 of Ref.~\cite{Fecko}.
    
    In this article we will use Gamma coefficient $\Gamma^c_{ab}$ and structure 
    coefficients $C^c_{ab}$.
    We recall that, for a given frame field $e_a$ and its related coframe field 
    $e^a$, the Gamma coefficients are defined by $\Gamma^c_{ab}=\langle 
    e^c,\nabla_be_a\rangle$, with the notation $\nabla_a = \nabla_{e_a}$ and 
    $\nabla$ the Levi-Civita connection.
    Then, by calculating $\nabla_b\langle e^c, e_a\rangle$, we find that 
    $(\nabla_b e^c)_a = -\Gamma^c_{ab}$. 
    The structure coefficients $C^c_{ab}$, appear through
    \begin{equation}
        [e_a,e_b] = C^c_{ab} e_c,\qquad
        d e^a = -\frac{1}{2}C^a_{bc} e^b\wedge e^c,\label{Eq. anholonomy}
    \end{equation}
    where we apply Einstein's summation convention for repeated indices, which 
    holds in the whole article unless explicitly stated otherwise.
    In the present article we need the following lemma.
    \begin{lemme}\label{Lem. C}
        In a given orthonormal frame $e_a$, with $e^a$ its coframe, and with 
        $\nabla$ the Levi-Civita connection we have
        \begin{equation*}
            C^a_{ab} = - \langle \nabla_a e^a, e_b\rangle = \Gamma^a_{ba}.
        \end{equation*}
    \end{lemme}
    \begin{proof}
        Indeed, for a given Levi-Civita connection $\nabla$, we have
        \begin{align*}
            g(\nabla_u v,w)
            = \tfrac{1}{2}\{&
            u(g(v,w)) + v(g(u,w)) - w(g(u,v))\\
            &+ g([u,v],w) - g([u,w],v) - g(u,[v,w])\},
        \end{align*}
        with $u$, $v$ and $w$ vectors.
        Now, with $u = v = e_a$ and $w = e_b$, with no summation, one gets
        \begin{equation*}
            g(\nabla_a e_a, e_b) = -\eta_{aa} C^a_{ab},\quad\text{(no 
            summation)}
        \end{equation*}
        in an orthonormal frame with $\nabla_{e_a} \equiv \nabla_a$ and 
        $\eta_{ab} = g(e_a,e_b) = \pm 1$ iff $a=b$ and zero otherwise.
        However,
        \begin{equation*}
            g(\nabla_a e_a, e_b)
            = \langle \flat \nabla_a e_a, e_b\rangle
            = \langle \nabla_a \flat e_a, e_b\rangle
            = \eta_{aa} \langle \nabla_a e^a, e_b\rangle,\quad\text{(no 
            summation)}.
        \end{equation*}
        Now, simplifying and allowing for the summation on $a$ yields the 
        result.
    \end{proof}
    
    In what follows the indices $l, m, p$ range from 1 to $n+k$ and refer to 
    the manifold $M$, the indices $\mu, \nu, \rho$ range from 1 to $n$ and 
    refer to the submanifold $S$ and, finally, the indices $A, B, C$ range from 
    $n+1$ to $n+k$ and refer to the normal bundle $S^\perp$.
    A frame $e_l$ of $M$ is said to be adapted to $S$ when, at every point of 
    $S$, the $e_\mu$ are in $T(S)$ and the $e_\sA$ in $S^\perp$.
    
    \section{Restriction of operators \label{Sec. Restr.}}
    We will now obtain the restriction of the codifferential by first 
    addressing the Hodge operator. 
    Special cases of the following formulas, with $k=2$ and using properties of 
    specific embedding formulas in an ambient null-cone, were found earlier in 
    a physics-motivated context \cite{Huguet:2025yll}.
    Here, these results reach their full generality.
    
    \subsection{Statements of the theorems}
    \begin{theorem}\label{Res_Hodge}
        Within the scope of the hypotheses \hyperlink{H}{(H)} (see 
        page~\pageref{(H)}), we have
        \begin{equation}
            \label{Eq. hodge}
            s^* *\alpha =(-1)^{nk} *_\sS i^{\omega^\sn}\alpha,
        \end{equation}
        at every point of $S$ and for all $\alpha\in\Omega(M)$, where $*$ and 
        $*_\sS$ are the Hodge operators on $M$ and $S$ respectively.
    \end{theorem}
    
    We deduce the following:
    \begin{theorem}\label{th : delta}
        Within the scope of the hypotheses \hyperlink{H}{(H)}, we have
        \begin{align}
            s^*\delta\alpha 
            &= \delta_\sS (s^*\alpha) 
            + \R(j^{\omega^\sn}\Theta_{\omega_\sn})\alpha \label{Eq. delta1}\\
            &= \delta_\sS (s^*\alpha)
            -i_{\omega_\sn}(\llbracket\omega^\sn,\alpha\rrbracket
            +(-1)^k (\delta\omega^\sn)\wedge\alpha),\label{Eq. delta}
        \end{align}
        at every point of $S$ and for all $\alpha\in\Omega(M)$, where $\delta$ 
        and $\delta_\sS$ are the codifferential operators on $M$ and $S$ 
        respectively.
    \end{theorem}

    The restriction of the Hodge Laplacian, defined by $\Delta = -(d\delta + 
    \delta d)$, follows easily:
    \begin{corollaire}
        Under the same hypotheses \hyperlink{H}{(H)}, we have
        \begin{multline}
            s^*\Delta\alpha
            = \Delta_\sS (s^*\alpha)
            +di_{\omega_\sn}(\llbracket\omega^\sn,\alpha\rrbracket+(-1)^k 
            \delta(\omega^\sn)\wedge\alpha) 
            +i_{\omega_\sn}(\llbracket\omega^\sn,d\alpha\rrbracket+(-1)^k 
            \delta(\omega^\sn)\wedge d\alpha),
            \label{Eq. square}
        \end{multline}
        which can be preferred in the form
        \begin{equation}
            s^*\Delta =\Delta_\sS s^*
            -\R(\delta j^{\omega^\sn}\Theta_{\omega_\sn}
            + j^{\omega^\sn}\Theta_{\omega_\sn}\delta).
            \label{Eq. square1}
        \end{equation}
    \end{corollaire}

    \begin{remark}
        The Hodge Laplacian is related to the rough Laplacian  
        $\nabla^*\nabla=\Tr(\nabla^2)$ thanks to the Weitzenböck formula (see 
        Ref.~\cite{Jost}) which, on differential forms, reads as
        \begin{equation}
            \label{eq.Weitzenbock}
            \nabla^*\nabla\alpha
            = \Delta\alpha + \mathrm{Rm}_{abcd}j^ai^bj^ci^d\alpha,
        \end{equation}
        with $\mathrm{Rm}$ the Riemann tensor,
        see the appendix~A.1. of Ref.~\cite{Huguet:2022rxi} and Eq.~(A7) of 
        Ref.~\cite{Huguet:2025yll}.
        Equations~\eqref{Eq. square} and \eqref{eq.Weitzenbock} can be used to 
        relate the ambient rough Laplacian on $M$ with the intrinsic rough 
        Laplacian on $S$.
    \end{remark}

    The compact form of these results raises the question of their practical 
    applicability.
    We show in Sec.~\ref{Sec. small a} that, when reformulated appropriately, 
    they enable concrete calculations, including the case $s^*\Delta\lambda$ 
    for $\lambda\in\Omega^1$.
    
    \subsection{Proofs of the theorems}
    
    In what follows, $e^\mu, \mu=1,\ldots,n$ is (locally) a direct orthonormal 
    frame of $T^*(S)$ and $e^\sA, A=n+1,\ldots,n+k$ is a direct orthonormal 
    basis of $S^\perp$.
    The orientations are chosen such that $e^l, l=1,\ldots,n+k$ is a direct 
    orthonormal frame of $T^*(M)$.
    For the sake of readability, we set $\eP = \sgn(g_\perp)$.
    Recall that, in this basis, we have $\omega^\sn=e^{n+1}\wedge\cdots\wedge 
    e^{n+k}$ and $\omega_\sn=e_{n+1}\wedge\cdots\wedge 
    e_{n+k}=\sgn(g^\perp)\sharp\omega^\sn = \eP\sharp\omega^\sn$ and 
    $i_{\omega_\sn} = i_{e_{n+k}}\cdots i_{e_{n+1}} \equiv i_{n+k}\cdots 
    i_{n+1}$ and $j^{\omega^\sn} = j^{e^{n+1}}\cdots j^{e^{n+k}} \equiv 
    j^{n+1}\cdots j^{n+k}$.
    
    Regarding the restriction to $S$, we have, at every point of $S$, 
    $s^*\phi=\phi$ for any scalar $\phi$ and, moreover, the identification 
    $s_*e_\mu=e_\mu$ allows the identification $s^* e^\mu=e^\mu$ even though 
    $s^*e^\sA=0$.
    As a consequence, of particular interest for the following we note, 
    developing $\alpha$ in the above basis, that
    $
    s^*\alpha 
    = i_{\omega_\sn} j^{\omega^\sn}\alpha
    = \eP i^{\omega^\sn}j^{\omega^\sn}\alpha$, and also 
    $s^*i_{\omega_\sn}=i_{\omega_\sn}$
    at every point of $S$.
    
    \begin{proof}[Proof of Theorem \ref{Res_Hodge}]
        First, note that for $\alpha_\sS$ a differential form spanned by $e^1, 
        \cdots, e^n$, we have
        \begin{align*}
            *j^{\omega^\sn}\alpha_\sS
            &= *j^{\omega^\sn\wedge\alpha_\sS}1
            = *\R(i^{\omega^\sn\wedge\alpha})1
            = **^{-1}i^{\alpha_\sS}i^{\omega^\sn}*1
            = i^{\alpha_\sS}i^{\omega^\sn}\omega_\sS\wedge \omega^\sn\\
            &=\eP (-1)^{nk}i^{\alpha_\sS}\omega_\sS
            =\eP (-1)^{nk}*_\sS\alpha_\sS,
        \end{align*}
        using the fact that $i^{\omega^\sn}\omega^\sn\wedge\alpha_\sS 
        = \eP i_{\omega_\sn}\omega^\sn\wedge\alpha_\sS
        = \eP\alpha_\sS$ and  
        using Eq.~(41) of Ref.~\cite{Huguet:2024mza}, namely $i^\alpha\omega = 
        *\alpha$ applied in $S$.
        Then, regarding the restriction of $*\alpha$ to $S$ we have 
        \begin{align*}
            s^**\alpha
            &= \eP i^{\omega^\sn}j^{\omega^\sn}*\alpha
            = \eP **^{-1}i^{\omega^\sn}**^{-1}j^{\omega^\sn}*\alpha\\
            &= \eP 
            *\R(i^{\omega^\sn})(-1)^{-k(a-k)}\R(j^{\omega^\sn})(-1)^{ka}\alpha\\
            &=\eP *j^{\omega^\sn}i^{\omega^\sn}\alpha\\
            &= (-1)^{nk}*_\sS i^{\omega^\sn}\alpha.\qedhere
        \end{align*}
    \end{proof}

    \begin{proof}[Proof of Theorem \ref{th : delta}]
        Note that the formula on the restriction of the Hodge operator 
        [Eq.~\eqref{Eq. hodge}] directly entails
        \begin{equation}
            \label{eq.Hodge-1}
            s^* *^{-1}\alpha
            = \eP (-1)^{k(a-1)} *_\sS^{-1} i^{\omega^\sn}\alpha.
        \end{equation}
        Now we have that
        \begin{align*}
            s^*\delta\alpha
            &= (-1)^a s^*(*^{-1}d*\alpha)
            = (-1)^a \tau *_\sS^{-1}i^{\omega^\sn}d*\alpha\\
            &= (-1)^a \tau *_\sS^{-1}(s^*\Theta^{\omega^\sn} + 
            (-1)^ks^*di^{\omega^\sn})*\alpha,
        \end{align*}
        using Eq.~\eqref{eq.Hodge-1} with the sign $\tau = \eP (-1)^{k(n+k-a)}$ 
        and using the fact that 
        $i^{\omega^\sn} d =s^*i^{\omega^\sn} d 
        = s^*[i^{\omega^\sn},d]^\g + (-1)^ks^*di^{\omega^\sn}
        = s^*\Theta^{\omega^\sn} + (-1)^ks^*di^{\omega^\sn}$.
        Regarding the first term we note that
        \begin{align*}
            \tau *_\sS^{-1}s^*\Theta^{\omega^\sn}*\alpha
            &= \eP \tau *_\sS^{-1}  
            i^{\omega^\sn}j^{\omega^\sn}\Theta^{\omega^\sn}*\alpha
            = \eP s^* *^{-1}j^{\omega^\sn}\Theta^{\omega^\sn}*\alpha\\
            &= \eP (-1)^a s^*\R(j^{\omega^\sn}\Theta^{\omega^\sn})\alpha,
        \end{align*}
        using Eq.~\eqref{eq.Hodge-1} in reverse.
        Regarding the second term we note that
        \begin{align*}
            *_\sS^{-1}s^*di^{\omega^\sn}*\alpha
            &= *_\sS^{-1}d_\sS s^***^{-1}i^{\omega^\sn}*\alpha
            = (-1)^{nk} *_\sS^{-1}d_\sS*_\sS  
            i^{\omega^\sn}(-1)^{ka}\R(i^{\omega^\sn})\alpha\\
            &= (-1)^{nk+ka+a}\delta_\sS  i^{\omega^\sn}j^{\omega^\sn}\alpha
            = \eP (-1)^{nk+ka+a} \delta_\sS s^*\alpha,
        \end{align*}
        using $s^*d=d_\sS s^*$.
        Thus, adding both terms and simplifying the signs, we get
        \begin{equation*}
            s^*\delta\alpha
            = \delta_\sS s^*\alpha
            + \eP s^*\R(j^{\omega^\sn}\Theta^{\omega^\sn})\alpha
            = \delta_\sS s^*\alpha
            - s^* i_{\omega_\sn}(\S^{\omega^\sn} + (-1)^k 
            j^{\delta\omega^\sn})\alpha,
        \end{equation*}
        which concludes the proof since $s^*i_{\omega_\sn}=i_{\omega_\sn}$.
    \end{proof}
    
    \section{Transverse differential forms}
    \label{Sec. transv}
    
    In this section, we put forward the notion of transverse and strongly 
    transverse differential forms.
    For these particular differential forms, the formulas obtained in 
    Sec.~\ref{Sec. Restr.} simplify significantly.
    
    \begin{definition}
        A differential form $\alpha\in\Omega(M)$ will be called transverse with 
        respect to $S$ if $i_v \alpha = 0$ for all $v\in S^\perp$.
        In addition, a differential form $\alpha\in\Omega(M)$ will be called 
        strongly transverse with respect to $S$ if $\alpha$ and $d\alpha$ are 
        transverse.
    \end{definition}
    
    \begin{remark}
        At every point of $S$ the differential form $\alpha$ is transverse iff 
        $s^*\alpha=\alpha$.
    \end{remark}
    \begin{remark}
        A strongly transverse form is, in some sense, an intrinsic form of the 
        submanifold $S$.
        In fact, in a coordinate system $y^m$ suitably adapted to the embedding 
        of $S$ in $M$ (i.e. $y\in S$ iff $y^\sA=0$), this implies that a 
        transverse form $\alpha$ reads as
        $\alpha=\alpha_{\mu_1\ldots\mu_a}dy^{\mu_1}\wedge\cdots\wedge 
        y^{\mu_a}/a!$, with no $dy^\sA$ in its expansion.
        In the case that it is strongly transverse we have in addition that 
        $\dr_\sA\alpha_{\mu_1\ldots\mu_a} = 0$, i.e. the components of $\alpha$ 
        do not depend on $y^\sA$. 
    \end{remark}
    
    Since $(s^*)^2 = s^*$, $s^*$ and $s^\perp=\Id-s^*$ together form a pair of 
    projections, we have $(s^\perp)^2 = s^\perp$ and $s^*s^\perp = s^\perp s^* 
    = 0$.
    Consequently, $\Omega(M) = \Im s^* \oplus \Im s^\perp$. Occasionally, we 
    will write $s^\perp\alpha = \alpha_\perp$ and $s^*\alpha=\alpha_\sS$. 
    $\Im s^*$ is the set of transverse forms and the forms belonging to $\Im 
    s^\perp$ are said to be longitudinal.
    
    We now introduce the vector field $v_s$ on $S$, which depends only on the 
    injection $s : S\to (M,g)$.
    This vector field will feature in the main formulas that follow.
    
    \begin{definition}
        Let 
        \begin{equation*}
            V_s(u,\lambda)=\sharp s^*(\nabla_u \lambda),\quad
            u\in T(S)^\perp,\quad
            \lambda\in T^*(S)^\perp,
        \end{equation*}
        which is tensorial in both entries and then is a tensor on $S^\perp$ 
        with values on $T(S)$.
        We now define
        \begin{equation}\label{eq:Ws}
            v_s = \Tr\ V_s = \sharp s^*(\nabla_\sA e^\sA)
            = g^{\mu\nu}(\nabla_\sA e^\sA)_\nu e_\mu 
            = -g^{\mu\nu}\Gamma^\sA_{\nu\sA}e_\mu
            = -\tilde{C}_\sA^{\sA\mu} e_\mu,
        \end{equation}
        as the trace of $V_s$.
        This is a vector field in $S$.
        Here, the second-to-last expression of $v_s$ holds in any adapted frame 
        and the
        last one holds in an orthonormal adapted frame with $\tilde{C}^{lm}_p = 
        \eta^{lq}\eta^{mr}\eta_{ps}C^s_{qr}$.
    \end{definition}

    \begin{property}\label{Prop. delta_transverse}
        Within the scope of the hypotheses \hyperlink{H}{(H)}, for $\alpha$ a 
        transverse differential form, at any point of $S$, we have
        \begin{equation}
            s^*\delta\alpha=\delta_\sS\alpha + i_{v_s}\alpha,\label{Eq. 
            delta_transv}
        \end{equation}
        recalling that at any point of $S$ we have $\alpha = s^*\alpha$ since 
        $\alpha$ is transverse.
    \end{property}
    \begin{proof}
        First, note that
        \begin{align*}
            s^*\delta
            &= s^*i_{\omega_\sn}j^{\omega^\sn}\delta
            = s^*\{\delta i_{\omega_\sn}j^{\omega^\sn}
            + [i_{\omega_\sn}j^{\omega^\sn},\delta]^\g\}\\
            &= s^*\{\delta i_{\omega_\sn}j^{\omega^\sn}
            - (-1)^k i_{\omega_\sn}[\delta, j^{\omega^\sn}]^\g 
            + (-1)^k\eP[i^{\omega^\sn},\delta]^\g j^{\omega^\sn}\}\\
            &= s^*\{\delta i_{\omega_\sn}j^{\omega^\sn}
            - i_{\omega_\sn}(\S^{\omega^\sn} + (-1)^k j^{\delta\omega^\sn})
            + (-1)^k\eP i^{d\omega^\sn}j^{\omega^\sn}\},
        \end{align*}
        using Eqs.~\eqref{eq:Graded_Bracket_1}--\eqref{eq:Graded_Bracket_3},
        the fact that $i_{\omega_\sn} = \eP i^{\omega^\sn}$ and
        the two brackets $[\delta, j^\alpha]^\g$ and $[i^\alpha,\delta]^\g$ 
        recalled in Table~\ref{tab:commutateurs}.
        Second, applying this identity to a transverse differential form 
        $\alpha$ and equating to Eq.~\eqref{Eq. delta} yields:
        \begin{equation}\label{eq:s^* delta a}
            s^*\delta\alpha
            = \delta_\sS (s^*\alpha)
            - (-1)^k\eP s^*i^{d\omega^\sn} j^{\omega^\sn}\alpha,
        \end{equation}
        since $i_{\omega_\sn}j^{\omega^\sn}\alpha = s^*\alpha = \alpha$ for a 
        transverse differential form.
        
        Now, let us assume that the frame $e_\sA$ is orthonormal, then
        \begin{align*}
            d\omega^\sn
            &= \sum_{i=1}^k (-1)^{i-1} e^{n+1}\wedge\cdots\wedge 
            de^{n+i}\wedge\cdots\wedge e^{n+k}\\
            &= \sum_{i=1}^k (-1)^{i-1} e^{n+1}\wedge\cdots\wedge
            \bigl(-C^{n+i}_{n+i\,\mu} e^{n+i}\wedge e^\mu
            - \tfrac{1}{2} C^{n+i}_{\mu\nu} e^\mu\wedge e^\nu
            \bigr )\wedge\cdots\wedge e^{n+k}\\
            &= (-1)^k\omega^\sn\wedge (C^\sA_{\sA \mu} e^\mu)
            - \tfrac{1}{2}C^\sA_{\mu\nu} e^\mu\wedge e^\nu\wedge i_\sA 
            \omega^\sn\\
            &= -(-1)^k\omega^\sn\wedge (\nabla_\sA e^\sA)
            - \tfrac{1}{2}C^\sA_{\mu\nu} e^\mu\wedge e^\nu\wedge i_\sA 
            \omega^\sn,
        \end{align*}
        since $C^\sA_{\sA \mu} = -(\nabla_{\sA} e^\sA)_\mu$, using 
        lemma~\ref{Lem. C}.
        
        Finally, with the above expression of $d\omega^\sn$, the last term in 
        the right-hand side of Eq.~\eqref{eq:s^* delta a} simplifies to
        \begin{align*}
            -(-1)^k \eP s^* i^{d\omega^\sn}j^{\omega^\sn}\alpha
            &= \eP s^* i^{\omega^\sn\wedge(\nabla_\sA 
            e^\sA)}j^{\omega^\sn}\alpha\\
            &= s^* i^{(\nabla_\sA e^\sA)}\eP i^{\omega^\sn}j^{\omega^\sn}\alpha
            = s^* i^{(\nabla_\sA e^\sA)}\alpha
            = i^{(\nabla_\sA e^\sA)}\alpha,
        \end{align*}
        as the term involving $e^\mu\wedge e^\nu$ necessarily vanishes since, 
        after computation, the term $j^\sA$ remains and is canceled by $s^*$, 
        and the last equality holds since $\alpha$ is transverse.
        Finally, with Eq.~\eqref{eq:Ws}, we recognize that $i^{(\nabla_\sA 
        e^\sA)}\alpha = i_{v_s}\alpha$ for $\alpha$ a transverse differential 
        form.
    \end{proof}

    \begin{remark}
        The above argument is a circuitous proof of Prop.~\ref{Prop. 
        delta_transverse}.
        A direct but tedious proof is also possible, using only the 
        Schouten–Nijenhuis bracket properties and Lemma~\ref{Lem:Lemma1}.
        The same is true for the proof of Theorem~\ref{Res_Hodge}, for which a 
        rather lengthy proof is available that depends only on the well-known 
        rules governing the action of the Hodge operator on orthonormal bases.
    \end{remark}

    \begin{remark}
        Property~\ref{Prop. delta_transverse} gives us the value of 
        $s^*\delta\alpha$ when $\alpha$ is transverse.
        Moreover, it turns out that $s^\perp\delta\alpha$ vanishes in this case.
        Therefore Eq.~\eqref{Eq. delta_transv} becomes
        \begin{equation*}
            \delta\alpha=\delta_\sS\alpha + i_{v_s}\alpha,%
        \end{equation*}
        for $\alpha$ transverse. This can be proved as follows.
        We start with
        \begin{equation*}
            s^*\delta
            = \delta s^*
            - i_{\omega_\sn}(\S^{\omega^\sn} + (-1)^k j^{\delta\omega^\sn})
            + (-1)^k\eP i^{d\omega^\sn}j^{\omega^\sn},
        \end{equation*}
        in which we assume that $\alpha$ is transverse and then project it onto 
        longitudinal forms, using $s^\perp=\Id-s^*$.
        We then obtain, in an orthonormal adapted frame and for $\alpha$ 
        transverse, that 
        \begin{equation*}\label{Eq. vertical}
            s^\perp (\delta\alpha)= \tfrac{1}{2} \tilde 
            C_\sA^{\mu\nu}e^\sA\wedge(i_\mu i_\nu \alpha),
        \end{equation*} 
        and, eventually, $C^\sA_{\mu\nu}=\langle e^\sA,[e_\mu,e_\nu]\rangle=0$ 
        since $[e_\mu,e_\nu]\in T(S)$.
    \end{remark}

    A direct consequence of Prop.~\ref{Prop. delta_transverse} is the following.
    \begin{corollaire}
        Let $\alpha$ be a strongly transverse differential form, from 
        Prop.~\ref{Prop. delta_transverse} we have
        \begin{equation*}
            s^*\Delta\alpha
            = \Delta_\sS \alpha
            - \L_{v_s}\alpha.
        \end{equation*}
    \end{corollaire}

    \section{Explicit expressions for low-degree differential forms} 
    \label{Sec. small a}
    
    In this section, we compute explicitly the restriction of the Hodge 
    Laplacian acting on scalar fields $\phi$ and  $1$-forms $\lambda$.
    To this end, we must evaluate $s^*\delta\alpha$ for $a = \deg\alpha = 1,2$ 
    since, writing
    \begin{equation*}
        s^* \delta\alpha
        =\delta_\sS(s^*\alpha) + T_a\alpha,\quad \alpha\in\Omega^a,
    \end{equation*}
    we have, for $\lambda\in\Omega^1$,
    \begin{equation}
        s^*\Delta\lambda =\Delta_\sS(s^*\lambda)- T_2(d\lambda)-d_\sS 
        T_1(\lambda),
        \label{Eq. Delta} 
    \end{equation}
    and for $\phi\in\Omega^0$
    \begin{equation}
        \Delta\phi =\Delta_\sS\phi - T_1(d\phi).
        \label{Eq. Delta phi} 
    \end{equation}
    Directly from Prop.~\ref{Prop. delta_transverse}, we find
    \begin{theorem}
        Within the scope of the hypotheses \hyperlink{H}{(H)}, for 
        $\alpha\in\Omega(M)$ and at every point of $S$, we have
        \begin{equation}
            s^*\delta\alpha
            = \delta_\sS(s^*\alpha) + i_{v_s}s^*\alpha
            + s^*\delta(\alpha_\perp),
            \label{Eq. small a}
        \end{equation}
        with $\alpha_\perp = s^\perp\alpha = \alpha - s^*\alpha$.
    \end{theorem}
    
    This theorem restates Prop.~\ref{Prop. delta_transverse}, but is useful 
    because, for small $a=\deg\alpha$, $s^*\delta(\alpha_\perp)$ can be made 
    explicit using an adapted frame.
    Indeed, in such a frame, one has, for $a=1$, $\lambda_\perp=\lambda_\sA 
    e^\sA$ and, for $a=2$,
    $\alpha_\perp = \alpha_{\mu\sA}e^\mu\wedge e^\sA
    + \frac{1}{2}\alpha_{\sA\sB}e^\sA\wedge e^\sB$.
    The codifferential of these terms can then be expressed in terms of the  
    components of $\lambda$ or $\alpha$ and, for instance, the coefficients 
    $\Gamma$ or, if the frame is orthonormal, the coefficients $C$.
    A direct computation yields
    \begin{align*}
        s^*\delta(\lambda_\perp)
        &= \lambda_\sA g^{lm}\Gamma^\sA_{lm} - g^{\sA\sB}e_\sA(\lambda_\sB), \\
        s^*\delta(\alpha_\perp)
        &= [\alpha_{\rho\sA}(g^{\rho\nu}\Gamma^\sA_{\mu\nu}
        - g^{\sA\sB}\Gamma^\rho_{\mu\sB})
        - \alpha_{\mu\sA}g^{lm}\Gamma^\sA_{lm}
        + g^{\sA\sB}e_\sB(\alpha_{\mu\sA})
        +  \alpha_{\sA\sB}g^{\sA\sC}\Gamma^\sB_{\mu\sC}]e^\mu,
    \end{align*}
    using the identity $\delta = -i^l\nabla_l$, the Leibniz formula for 
    $\nabla_l$, $s^*e^\mu = e^\mu$, $s^*e^\sA = 0$, $s^* g^{\mu\sA} = 0$ and 
    $\nabla_l e^p = -\Gamma^p_{ml}e^m$.
    In   an orthonormal basis it simplifies to
    \begin{align*}
        s^*\delta(\lambda_\perp) 
        &= -\lambda_\sA \tilde{C}^{p\sA}_p
        - \eta^{\sA\sB}e_\sA(\lambda_\sB),\\
        s^*\delta(\alpha_\perp)
        &= [\alpha_{\mu\sA}\tilde{C}^{p\sA}_p
        - \alpha_{\nu\sA}\widetilde{C}^{\nu\sA}_\mu
        + \eta^{\sA\sB} e_{\sB}(\alpha_{\mu\sA})
        - \tfrac{1}{2}\alpha_{\sA\sB}\widetilde{C}^{\sA\sB}_\mu] e^\mu,
    \end{align*}
    with $C^l_{pm}$ the anholonomy coefficients given in Eq.~\eqref{Eq. 
    anholonomy} and $\tilde{C}^{lm}_p = \eta^{lq}\eta^{mr}\eta_{ps}C^s_{qr}$, 
    using Eq.~\eqref{delta_j},
    $\llbra e^l, e^m\rrket = \flat\L_{\sharp e^l} \sharp e^m = \flat [\sharp 
    e^l, \sharp e^m] = \tilde{C}^{lm}_p e^p$, 
    $\llbracket\phi,\beta\rrbracket=-i^{d\phi}\beta$ (see 
    Ref.~\cite{Huguet:2024mza}), and  $\delta e^l = -\eta^{lm}C^p_{pm}$ (see 
    App.~A.2 of Ref.~\cite{Huguet:2022rxi}).
    This yields, using Eqs.~\eqref{eq:Ws} and \eqref{Eq. small a}, the 
    following theorem.
    
    \begin{theorem}\label{Prop. delta small a}
        Within the scope of the hypotheses \hyperlink{H}{(H)}, for 
        $\lambda\in\Omega^1(M)$, $\alpha\in\Omega^2(M)$ and at every point of 
        $S$, we have 
        \begin{align}
            s^*\delta\lambda
            &= \delta_\sS(s^*\lambda)
            - g^{\mu\nu}\Gamma^\sA_{\mu\sA}\lambda_\nu
            + \lambda_\sA g^{lm}\Gamma^\sA_{lm}
            - g^{\sA\sB}e_\sA(\lambda_\sB) \label{eq:First}\\
            &= \delta_\sS(s^*\lambda)
            - \tilde C^{\sA\mu}_\sA\lambda_\mu
            - \tilde C_p^{p\sA}\lambda_\sA
            - \eta^{\sA\sB}e_\sA(\lambda_\sB),\\
            s^*\delta\alpha
            & =\delta_\sS(s^*\alpha) \nonumber\\ &\quad
            + [\alpha_{\mu\nu}\Gamma^\sA_{\rho\sA}g^{\rho\nu}
            + \alpha_{\rho\sA}(g^{\rho\nu}\Gamma^\sA_{\mu\nu}
            - g^{\sA\sB}\Gamma^\rho_{\mu\sB}
            - \delta^\rho_{\mu}g^{lm}\Gamma^\sA_{lm})
            + g^{\sA\sB}e_\sB(\alpha_{\mu\sA})
            + \alpha_{\sA\sB}g^{\sA\sC}\Gamma^\sB_{\mu\sC}]e^\mu\\
            &= \delta_\sS(s^*\alpha)
            + [\alpha_{\mu\nu}\tilde C^{\sA\nu}_\sA 
            + \alpha_{\mu\sA}\tilde C_p^{p\sA}      
            - \alpha_{\nu\sA}\widetilde{C}^{\nu\sA}_\mu
            + \eta^{\sA\sB} e_{\sB}(\alpha_{\mu\sA})
            - \tfrac{1}{2}\alpha_{\sA\sB}\widetilde{C}^{\sA\sB}_\mu] e^\mu,
        \end{align}
        in which the first lines hold in an arbitrary adapted frame while  the 
        second ones  hold in an orthonormal adapted frame.
    \end{theorem}
    This property, together with Eqs.~\eqref{Eq. Delta} and \eqref{Eq. Delta 
    phi}, yields
    \begin{corollaire}
        \label{Prop. Delta small a}
        Within the scope of the hypotheses \hyperlink{H}{(H)}, for 
        $\phi\in\Omega^0(M)$, $\lambda\in\Omega^1(M)$, and at every point of 
        $S$, we have, in any adapted frame
        \begin{align}
            \Delta\phi&=\Delta_\sS\phi+g^{\mu\nu}\Gamma^\sA_{\mu\sA}e_\nu(\phi)-g^{lm}\Gamma^\sA_{lm}
             e_\sA(\phi)+ g^{\sA\sB}e_\sA(e_\sB(\phi)),\\
            s^*\Delta\lambda
            &=\Delta_\sS(s^*\lambda)
            + e_\mu(g^{\rho\nu}\Gamma^\sA_{\rho\sA}\lambda_\nu
            - \lambda_\sA g^{lm}\Gamma^\sA_{lm}
            + g^{\sA\sB}e_\sA(\lambda_\sB))e^\mu
            - (d\lambda)_{\mu\nu}\Gamma^\sA_{\rho\sA}g^{\rho\nu}e^\mu
            \nonumber\\
            &\quad -[(d\lambda)_{\rho\sA}(g^{\rho\nu}\Gamma^\sA_{\mu\nu}
            - g^{\sA\sB}\Gamma^\rho_{\mu\sB}
            - \delta^\rho_{\mu}g^{lm}\Gamma^\sA_{lm})
            + g^{\sA\sB}e_\sB((d\lambda)_{\mu\sA})
            + (d\lambda)_{\sA\sB}g^{\sA\sC}\Gamma^\sB_{\mu\sC}]e^\mu,
        \end{align}
        where 
        $(d\lambda)_{ml}=e_m(\lambda_l)-e_l(\lambda_m)-\lambda_p(\Gamma^p_{lm}-\Gamma^p_{ml})$.
        While, if the frame is, in addition, orthonormal, the equations become
        \begin{align}
            \Delta\phi
            &= \Delta_\sS\phi
            + \tilde C^{\sA\mu}_\sA e_\mu(\phi)
            + \tilde C_p^{p\sA}e_\sA(\phi)+ \eta^{\sA\sB}e_\sA(e_\sB(\phi)),\\
            s^*\Delta\lambda
            &= \Delta_\sS(s^*\lambda)
            + e_\mu(\tilde C^{\sA\nu}_\sA\lambda_\nu
            + \lambda_\sA \tilde C_p^{p\sA} + 
            \eta^{\sA\sB}e_\sA(\lambda_\sB))e^\mu
            \nonumber\\
            &\quad
            - [(d\lambda)_{\mu\nu}\tilde C^{\sA\nu}_\sA 
            +(d\lambda)_{\mu\sA}\tilde C_p^{p\sA}
            - (d\lambda)_{\nu\sA}\widetilde{C}^{\nu\sA}_\mu
            + \eta^{\sA\sB} e_{\sB}((d\lambda)_{\mu\sA})
            - \tfrac{1}{2}(d\lambda)_{\sA\sB}\widetilde{C}^{\sA\sB}_\mu] e^\mu,
            \label{eq:Last}
        \end{align}
        where $(d\lambda)_{ml}=e_m(\lambda_l)-e_l(\lambda_m)+\lambda_p 
        C^p_{ml}$.
    \end{corollaire}

\end{document}